\documentclass[12pt]{amsart}
\usepackage{stmaryrd}
\usepackage{mathrsfs}
\usepackage{amssymb}

\usepackage{titletoc}
\input xy
  \xyoption{all}
\usepackage{amscd}
\usepackage{amsmath, amssymb}
\usepackage{amsfonts}
\usepackage[colorlinks,linkcolor=blue,citecolor=blue, pdfstartview=FitH]{hyperref}

\numberwithin{equation}{section}

\theoremstyle{plain}
\newtheorem{thm}{Theorem}[section]
\newtheorem{theorem}[thm]{Theorem}
\newtheorem{lemma}[thm]{Lemma}
\newtheorem{corollary}[thm]{Corollary}
\newtheorem{proposition}[thm]{Proposition}

\theoremstyle{definition}

\newtheorem{remark}[thm]{Remark}

\newtheorem{definition}[thm]{Definition}

\newtheorem{example}[thm]{Example}

\newtheorem{defn-thm}[thm]{Definition-Theorem}

\newcommand{\sJ}{{\mathcal J}}

\newcommand{\sL}{{\mathcal L}}

\newcommand{\sO}{{\mathcal O}}

\newcommand{\C}{{\mathbb C}}

\renewcommand{\P}{{\mathbb P}}
\newcommand{\Q}{{\mathbb Q}}

\newcommand{\Z}{{\mathbb Z}}

\newcommand{\rank}{{rank}}

\newcommand{\Ker}{{ Ker}}

\newcommand{\Ext}{{ Ext}}

\newcommand{\ds}{\oplus}
\newcommand{\bds}{\bigoplus}
\newcommand{\ts}{\otimes}

\newcommand{\btheorem}{\begin{theorem}}
\newcommand{\etheorem}{\end{theorem}}
\newcommand{\bproposition}{\begin{proposition}}
\newcommand{\eproposition}{\end{proposition}}
\newcommand{\bdefinition}{\begin{definition}}
\newcommand{\edefinition}{\end{definition}}
\newcommand{\bcorollary}{\begin{corollary}}
\newcommand{\ecorollary}{\end{corollary}}
\newcommand{\bproof}{\begin{proof}}
\newcommand{\eproof}{\end{proof}}
\newcommand{\bremark}{\begin{remark}}
\newcommand{\eremark}{\end{remark}}
\newcommand{\eexample}{\end{example}}
\newcommand{\bexample}{\begin{example}}

\newcommand{\elemma}{\end{lemma}}
\newcommand{\blemma}{\begin{lemma}}

\renewcommand{\bar}{\overline}
\newcommand{\eps}{\varepsilon}

\renewcommand{\phi}{\varphi}

\newcommand{\ee}{\end{eqnarray*}}
\newcommand{\be}{\begin{eqnarray*}}

\newcommand{\beq}{\begin{equation}}
\newcommand{\eeq}{\end{equation}}

\newcommand{\bd}{\begin{enumerate}}
\newcommand{\ed}{\end{enumerate}}

\renewcommand{\tilde}{\widetilde}

\renewcommand{\>}{\rightarrow}

\usepackage{fancyhdr}
\renewcommand{\leftmark}{{\scriptsize Characterizations of projective spaces and quadrics}}

\begin{document}
\title{Characterizations of projective spaces and quadrics by strictly nef bundles}
\makeatletter
\let\uppercasenonmath\@gobble
\let\MakeUppercase\relax
\let\scshape\relax
\makeatother
\author{ Duo Li and Xiaokui Yang}

\address{{Address of Duo Li: Yau Mathematical Sciences Center, Tsinghua University, Beijing, 100084, P. R. China.}}
\email{\href{mailto:duoli@math.tsinghua.edu.cn}{{duoli@math.tsinghua.edu.cn}}}
\address{{Address of Xiaokui Yang: Morningside Center of Mathematics, Institute of
        Mathematics, Hua Loo-Keng Key Laboratory of Mathematics,
        Academy of Mathematics and Systems Science,
        Chinese Academy of Sciences, Beijing, 100190, China.}}
\email{\href{mailto:xkyang@amss.ac.cn}{{xkyang@amss.ac.cn}}}
\maketitle

\begin{abstract} In this paper,  we show that if the
tangent bundle of a smooth projective variety is strictly nef, then
it is isomorphic to a projective space; if a projective variety
$X^n$ $(n>4)$ has  strictly nef $\Lambda^2 TX$, then it is
isomorphic to $\P^n$ or quadric $\Q^n$. We also prove that on
elliptic curves, strictly nef vector bundles are ample, whereas
there exist Hermitian flat and strictly nef vector bundles on any
smooth curve with genus $g\geq 2$.

\end{abstract}

\setcounter{tocdepth}{1} \tableofcontents

\section{Introduction}

\noindent Let $X$ be a smooth projective variety over $\C$. A line
bundle $L$ is said to be strictly nef if $L\cdot C>0$ for any
irreducible curve $C\subset X$. A vector bundle $V$ is said to be
strictly nef (resp. ample, nef) if the tautological line bundle
$\sO_V(1)$ of the projective bundle $\P(V)\>X$ is strictly nef
(resp. ample, nef). It is obvious that the notion of strictly
nefness is stronger than that of nefness, but weaker than that of
ampleness. It is well-known that  there are equivalent analytical
and cohomological characterizations for ampleness and nefness (e.g.
\cite{Dem}, \cite{Har66}) which play fundamental roles in complex
analytic geometry and algebraic geometry. However, there are no such
characterizations for strictly nef bundles.

  The purpose of this paper
is to investigate  properties of strictly nef vector bundles. In
particular, we want to show the similarity and differences between
various positivity notions. It is obvious that numerically trivial
line bundles are nef but not strictly nef. As a striking example, we
show there exist \emph{Hermitian flat} vector bundles which are
strictly nef (Example \ref{example}), which  is significantly
different from our primary impression on strictly nef vector
bundles. However, there are still some similarities between strictly
nefness and ampleness.

\btheorem\label{main1} Let $X$ be a smooth projective variety. If
$\Lambda^{r}TX$ is strictly nef, then $X$ is uniruled. \etheorem

\noindent In particular, if $\Lambda^r TX$ is strictly nef, then
$-K_X$ can not be numerically trivial. This is quite different from
the Hermitian flat and strictly nef vector bundle in Example
\ref{example}. Moreover, the existence of rational curves on such
varieties is  one of the key ingredients and also the starting point
in the characterizations of $\P^n$ or $\Q^n$ by strictly nef
bundles, since as a priori, we do not know whether the strictly
nefness of $\Lambda^rTX$ could imply the ampleness of $-K_X$.
Indeed,  a well-known conjecture of Campana and  Peternell
(\cite[Problem ~11.4]{CP91}) predicts
 that if $-K_X$ is strictly nef, then $-K_X$ is ample. This conjecture is only verified by Maeda for surfaces
(\cite{Mea93}) and by Serrano  for threefolds (\cite{Ser95}) (see
also \cite{Ueh00} and \cite{CCP06}). The proof of Theorem
\ref{main1} relies on a metric version (\cite{HW12}) of the
fundamental work (\cite{BDPP13}) of
Boucksom-Demailly-Peternell-Paun, Yau's solution to the Calabi
conjecture (\cite{Yau78}) and Yau's  criterion for Calabi-Yau
manifolds (e.g.\cite[Theorem~4]{Yau77}).

 As analogous
to Mori's fundamental result (\cite{Mor79}), we prove
 \btheorem\label{main2} Let
$X$ be a smooth projective variety with dimension $n$. If $TX$ is
 strictly nef, then  $X\cong \P^n$. \etheorem

\noindent That means, as holomorphic tangent bundle,  $TX$ is
strictly nef if and only if $TX$ is ample. Along the same line as
Theorem \ref{main2}, we obtain the following classification (see
also the classical results \cite{CS95} of Cho-Sato, \cite{Miy04}  of
Miyaoka and the recent work \cite{DH} of Dedieu-Hoering).

\btheorem\label{main3} Let $X$ be a smooth projective variety of
complex dimension $n>4$. Suppose $\Lambda^2TX$ is strictly nef, then
$X$ is isomorphic to $\P^n$ or quadric $\Q^n$.
 \etheorem

\noindent Note that, as a priori, the strictly nefness of $TX$ may
not imply the strictly nefness of $\Lambda^{2}TX$ (c.f. the
Hermitian flat and strictly nef rank $2$ bundle in Example
\ref{example}). This is also the key difficulty in the
characterizations of $\P^n$ or $\Q^n$ by strictly nef bundles. There
are also some other classical characterizations of $\P^n$ or $\Q^n$,
e.g. \cite{KO73}, \cite{SY80}, \cite{Siu80}, \cite{Mok88},
\cite{YZ90}, \cite{Pet90, Pet91}, \cite{AW01}, \cite{CMSB02},
\cite{ADK08}, \cite{Hwa13}, \cite{Yan16} and etc. We refer to the
papers \cite{Pet96, DH}, and the references therein.\vskip
0.4\baselineskip

 In Example \ref{example}, we show that, on any
smooth curve with genus $\geq 2$, there exist \emph{Hermitian flat}
vector bundles which are also strictly nef. However, we prove that
it can not happen on elliptic curves. Moreover, \btheorem
\label{main5} Let $V$ be a vector bundle over an elliptic curve $C$.
If  $V$ is strictly nef, then $V$ is ample. \etheorem \noindent Note
that Theorem \ref{main5} is also essentially used in the proof of
Theorem \ref{main3}.

 \vskip 0.3\baselineskip

 \noindent\textbf{Remark.} It is
not hard to see that the strictly nef condition in the above
theorems can be slightly relaxed  by using refined techniques
developed in \cite{BDPP13}, \cite{Hwa00} and \cite{Mok08}.

\vskip 0.3\baselineskip
 \noindent \textbf{Sketched ideas in the
proofs.} Since strictly nef bundles have no analytical or
coholomogical characterizations, the proofs of our theorems are
significantly different from those classical results. For instance,
from  Example \ref{example},  as a priori, the strictly nefness
$\Lambda^{2}TX$ can not imply the strictly nefness of $-K_X$,
nevertheless the ampleness of $-K_X$. Hence, we can not use
birational techniques on Fano manifolds. In the proof of Theorem
\ref{main3}, the key difficulty is to prove $-K_X$ is ample. Here we
argue by contradiction. Suppose $\Lambda^2TX$ is strictly nef and
$c_1^n(X)=0$. We divided it into two cases. (1) When $TX$ is not
nef. We show there exists an extremal rational curve (Theorem
\ref{main1}) such that $-K_X|_C$ is ample. Mover, by the ``Cone
Theorem" and the ideas in \cite[Theorem ~4.2]{SW}, we show the Mori
elementary contraction $\beta:X\longrightarrow B$ is indeed a fiber
bundle over a curve (Lemma \ref{bend and break} and Lemma
\ref{smooth fibration}). By analyzing  vector bundles over curves
carefully, we get a contradiction (Theorem \ref{nef}).  (2) When
$TX$ is nef. We use the structure theorem of
Demailly-Peternell-Schneider \cite[Main Theorem]{DPS94} and show
that $X$ has a finite \'etale cover of the form of a projective
bundle over an elliptic curve. Based on the key fact that strictly
nef bundles on elliptic curves are ample (Theorem \ref{main5}), we
are able to obtain a contradiction (Theorem \ref{keytheorem}),
thanks to a Barton-Kleiman type criterion for strictly nef bundles.

\vskip 0.3\baselineskip

\textbf{Acknowledgements.}  The
 authors  would like to thank Professor S.-T. Yau
for his comments and suggestions on an earlier version of this paper
which clarify and improve the presentations. The authors would also
like to thank Yifei Chen, Yi Gu  and   Xiaotao Sun  for some useful
discussions. The first author is very grateful to Professor Baohua
Fu for his support, encouragement and stimulating discussions over
the last few years.
   The second author
wishes to
 thank  Kefeng Liu, Valentino Tosatti and Xiangyu Zhou
 for  helpful discussions.
This work was partially supported by China's Recruitment
 Program of Global Experts and National Center for Mathematics and Interdisciplinary Sciences,
 Chinese Academy of Sciences.

\vskip 2\baselineskip

\section{Basic properties of   strictly nef vector bundles}

Let $X$ be a smooth projective variety. There are many equivalent
definitions for ampleness of  a line bundle (e.g. \cite{Har66},
\cite{Laz04I}). The Nakai-Moishezon-Kleiman criterion asserts  that:
a line bundle $L$ is \emph{ample} if and only if \beq L^{\dim
Y}\cdot Y>0 \eeq for every positive-dimensional irreducible
subvariety $Y\subset X$. Similarly, a line bundle $L$ is \emph{nef}
if and only if $ L\cdot C\geq 0 $ for every
 irreducible curve $C\subset X$. As an analogue,

\bdefinition   A line bundle $L$ is said to be \emph{strictly nef},
if for any irreducible curve $C$ in $X$, \beq L\cdot C>0.\eeq
\edefinition

\noindent There are many strictly nef line bundles which are not
ample (e.g.  Example \ref{example}).

  A vector bundle $V$ is said to be \emph{ample} (resp. \emph{nef}) if
  the tautological line bundle $\sO_V(1)$ of $\P(V)\>X$ is   ample
  (resp. nef). Similarly,

\bdefinition  A vector bundle $V$ is called \emph{strictly nef}, if
its tautological line bundle $\sO_V(1)$ of $\P(V)\>X$ is a strictly
nef line bundle. \edefinition

\noindent We have the following characterization of strictly nef
vector bundles which is analogous to the Barton-Kleiman criterion
for nef vector bundles (e.g. \cite[Proposition~6.1.18]{Laz04II}):

\blemma\label{cri0} Let $X$ be a smooth projective variety and $V$ a
holomorphic vector bundle. Then the following are equivalent: \bd
\item $V$ is  strictly nef;

\item for any finite morphism $\nu:C\>X$ where $C$ is a smooth curve, and any line bundle quotient $\nu^*(V)\twoheadrightarrow L$, one has
\beq \deg_C(L)>0.\label{wsn}\eeq \ed

\elemma

\bproof  It is well-known that giving a line bundle quotient
$\nu^*V\twoheadrightarrow L$ is the same as giving a map
$\mu:C\>\P(V)$ commuting with the projection to $X$: \beq \xymatrix{
  C \ar[rr]^{\mu} \ar[dr]_{\nu}
                &  &    \P(V)\ar[dl]^{\pi}    \\
                & X                }
\label{diag}\eeq with $L=\mu^*\left(\sO_V(1)\right)$, where
$\sO_V(1)$ is the tautological line bundle.

$(2)\Longrightarrow (1).$
 Let $C$ be
an arbitrary irreducible curve in $Y=\P(V)$. Let $\mu:C\>Y$ be the
inclusion and $\nu:C\>X$ the map given in (\ref{diag}). Then $$
\sO_V(1)\cdot C=\deg_C(L)>0. $$

$(1)\Longrightarrow (2).$ Let $\nu:C\>X$ be a finite morphism where
$C$ is an irreducible curve.  Let $\mu:C\>\P(V)$ be the map in
(\ref{diag}) and $\tilde C= \mu(C)$. Then \be \sO_V(1)\cdot \tilde
C&=&\deg_{ C} \left(\mu^*(\sO_V(1))\right)=\deg_C(L)>0\ee since
$\sO_V(1)$ is strictly nef. \eproof

\subsection{An  example of Hermitian flat and strictly nef vector bundle} In this subsection, we will show by an
example (essentially constructed  by D. Mumford, e.g.
\cite[Section~10, Chapter I]{Har70}) that the terminologies
introduced above are mutually different.

\bexample\label{example}  Let $C$ be a smooth curve of genus $g\geq
2$. There exists a  stable vector bundle $V$ such that $V$ is
\emph{Hermitian flat} and \emph{ strictly nef}. In particular, the
tautological line bundle $\sO_V(1)$ of the projective bundle
$\P(V)\>C$ is strictly nef, but it is not ample. \eexample

\noindent At first, we need to introduce the concept of
(semi-)stability for vector bundles.

\bdefinition A holomorphic  vector bundle $V$ over a smooth curve
$C$ is called stable(resp. semi-stable) if for any subbundle $V_1$
of $V$ with $0<\rank(V_1)<\rank(V)$, we have \beq
\frac{\deg(V_1)}{\rank (V_1)}<\frac{\deg(V)}{\rank (V)}\ \
\left(\text{resp.}\ \ \frac{\deg(V_1)}{\rank (V_1)}\leq
\frac{\deg(V)}{\rank (V)}\right).\eeq \edefinition

\blemma\label{Ses} Let $C$ be a smooth curve of genus $g\geq 2$.
Then there exists a rank two stable vector bundle $V$ of degree zero
such that all its symmetric powers $Sym^m V$ are stable. \elemma

\bproof This result is due to C.S. Seshadri. We refer the proof to
\cite[Theorem ~10.5, Chapter~I]{Har70}. \eproof

\bcorollary  The vector bundle $V$ in Lemma \ref{Ses} is Hermitian
flat. \ecorollary

\bproof The proof is simple and we include it here for readers'
convenience. By Hartshorne's Theorem (e.g.
\cite[Theorem~6.4.15]{Laz04II}), a stable vector bundle $V$ over a
curve is nef if and only if $\deg(V)\geq 0$. Since $V$ is stable and
of degree zero, $V$ is nef. On the other hand, since $V$ is of rank
$2$,  $V=V^*\ts \det V$. Hence $V^*\ts \det V$ is a nef vector
bundle. Since $\det V$ is numerically trivial, i.e. $\deg(V)=0$, we
see $V^*$ is also nef. Since $V$ is stable over $C$, there exists a
Hermitian-Einstein metric $h$ on $V$ (e.g. \cite{UY86}), i.e. \beq
g^{-1}\cdot R_{\alpha\bar \beta}= c\cdot h_{\alpha\bar \beta}\eeq
for some constant $c$ where $g$ is a smooth metric on  $C$. Since
both $V$ and $V^*$ are nef, we deduce $c=0$, i.e. $V$ is Hermitian
flat.
 \eproof

 \noindent Let $X=\P(V)$ and $\pi:\P(V)\>C$ be the projection. Let
 $\sO_V(1)$ be the tautological line bundle of $\P(V)$ and $D$ be
 the corresponding divisor over $X$. For any effective curve $Y$ on
 $X$, we denote by $m(Y)$ the degree of $Y$ over $C$. Then there
 is an exact sequence
 \beq   0\> Pic(C)\stackrel{\pi}{\>}Pic(X)\stackrel{m}{\>}\Z\>0.\eeq
It follows that  the divisors on $X$, modulo numerical equivalence,
form a free abelian group of rank $2$, generated by $D$ and $F$
where $F$ is any fiber of $\P(V)$. It is easy to see that \beq
D^2=\deg(V)=0,\ \ D\cdot F=1,\ \   F^2=0. \eeq

\blemma\label{examplelemma} For any $m>0$, there is a $1-1$
correspondence between \bd
\item effective curves $Y$ on $X$, having no fibers as components,
of degree $m$ over $C$; and

\item line subbundles $L$ of $Sym^m V$. \ed Moreover, under this
correspondence, one has \beq D\cdot Y=
m\deg(V)-\deg(L).\label{examplekey}\eeq \elemma

\bproof See \cite[Proposition~10.2, Chapter~I]{Har70}. \eproof

\noindent Note that
 $Sym^{ m} V$ are stable and of degree zero for all $m\geq 1$. For a line subbundle $L$ of $Sym^m V$, we
have
$$\deg(L)<\frac{\deg\left(Sym^mV\right)}{\rank\left(Sym^mV\right)}=0.$$

\noindent Let $Y$ be an arbitrary irreducible curve on $X$. If $Y$
is a fiber, then $D\cdot Y=1$. If $Y$ is an irreducible curve of
degree $m>0$ over $C$, then by Lemma \ref{examplelemma}, $Y$ is
corresponding to a line  subbundle $L$ of $Sym^mV$. Therefore, by
formula (\ref{examplekey}) \beq D\cdot
Y=m\deg(V)-\deg(L)=-\deg(L)>0.\eeq Hence, the line bundle $\sO_V(1)$
of the divisor $D$ is strictly nef, i.e. $V$ is  strictly nef. Since
$V$ is  flat, $V$ can not be ample. In particular, $\sO_V(1)$ is
strictly nef but not ample. \qed

\subsection{Basic properties for  strictly nef vector bundles} As we
have seen from the previous example, the terminology ``strictly nef"
is significantly different from  other notions. Many functorial
properties do not hold for strictly nef bundles. For instance, if
$E$ is strictly nef, $\Lambda^s E$ and $\det E$ are not necessarily
strictly nef. In Example \ref{example}, the strictly nef vector
bundle is flat, and so its determinant line bundle is flat and it is
not strictly nef. In this subsection, we prove several basic
properties for strictly nef vector bundles.

\bproposition\label{cricurve} $E$ is  strictly nef if and only if
$E|_C$ is strictly nef for every curve $C$. \eproposition \bproof It
follows from Lemma \ref{cri0}. \eproof

\bproposition  If $Sym^kE$ is  strictly nef for some $k\geq 1$, then
$E$ is  strictly nef.  \eproposition \bproof Let
$v_k:\P(E)\>\P(Sym^kE)$ be the Veronese embedding. Then we have \beq
\sO_E(k)=v_k^*(\sO_{Sym^kE}(1))\eeq Therefore  $E$ is  strictly nef
if $Sym^kE$ is  strictly nef for some $k\geq 1$.\eproof

\bproposition\label{directsum} \bd \item If $E$ is  strictly nef,
then any quotient of it is  strictly nef.

\item  If $E\ds F$  is  strictly nef, then both $E$ and $F$ are
 strictly nef.

\item  If $E$ is  strictly nef,   for any $L\in Pic^0(X)$, $E\ts
L$ is  strictly nef. \ed\eproposition

\bproof $(1)$. Let $Q$ be a quotient of $E$. Suppose $Q$ is not
strictly nef, by Lemma \ref{cri0}, there exist a  finite morphism
$\nu:C\>X$ where $C$ is a smooth curve, and a line bundle quotient
$\nu^*(Q)\twoheadrightarrow L$ such that  $\deg(L)\leq 0$. Note that
$L$ is also a quotient of $\nu^*(E)$, which is a contradiction.
$(2)$ follows from $(1)$.  $(3)$. Since $\P(E\ts L)= \P(E)$, \beq
\sO_{E\ts L}(1)=\sO_E(1)\ts \pi^*(L).\eeq Hence, if $L\in Pic^0(X)$,
$\sO_{E\ts L}(1)$ is strictly nef if and only if $\sO_E(1)$ is
strictly nef.
 \eproof

\begin{proposition}Let $X$, $Y$ be two smooth projective varieties and
$f:X\>Y$ be a  surjective morphism. If $f^*E$ is strictly nef, then
$E$ is  strictly nef.
\end{proposition}

\bproof Let $C$ be any curve in $Y$ and $\tilde C$ be a curve in $X$
such that $f(\tilde C)=C$. Since $f^*E|_{\tilde C}$ is strictly nef,
 $E|_C$ is strictly nef. By Proposition \ref{cricurve}, $E$ is
 strictly nef.
\eproof

\bcorollary\label{a} If $V$ is a  strictly nef bundle over $X$, then
$H^0(X,V^*\ts L)=0$ for any $L\in Pic^0(X)$. \ecorollary

\bproof  Let $\tilde V=V\ts L^*$. Since $L\in Pic^0(X)$ and $V$ is
 strictly nef, $\tilde V$ is  strictly nef. By
\cite[Proposition~1.16]{DPS94}, if $\tilde V$ is nef and
$H^0(X,\tilde V^*)\neq 0$,  there exists a nowhere vanishing section
$\sigma\in H^0(X,\tilde V^*)$. There is  a bundle map induced by
$\sigma$,
$$f_\sigma: X\times \C\>\tilde V^*,\ \ \ f_\sigma(x,v)=v\cdot\sigma(x).$$
In particular, the trivial bundle $\sO=X\times \C$ is a holomorphic
subbundle of $\tilde V^*$. By duality, the trivial bundle $\sO$ is a
quotient of $\tilde V$. By Lemma \ref{cri0}, $\tilde V$ can not be
 strictly nef. This is a contradiction. \eproof

\vskip 1\baselineskip

\section{Strictly nef bundles over elliptic curves}

\subsection{Projective bundles over Calabi-Yau manifolds}

 We need another characterization
of strictly nef anti-canonical bundles which is analogues to Lemma
\ref{cri0}.

\blemma \label{cri} Let $X=\P(V)$ be a projective bundle of a
holomorphic vector bundle $V$ over a Calabi-Yau manifold $Y$, i.e.
$K_Y\cong\sO_Y$. Then the following are equivalent:

\bd \item $-K_X$ is  strictly nef;

\item for any finite morphism $\nu:C\>Y$ from an irreducible curve $C$ to $Y$ and any line bundle quotient $\nu^*(V)\twoheadrightarrow L$, one has
\beq \deg(L)> \frac{\deg(\nu^*(V))}{\text{rank}(V)}.\label{sta}\eeq
\ed \elemma

\bproof Here we use the same configuration as in Lemma \ref{cri0}.

\noindent $(2)\Longrightarrow (1).$ Let $C$ be  an arbitrary
irreducible curve in $X=\P(V)$. Let $\mu:C\>X$ be the inclusion and
$\nu:C\>Y$ the map given in (\ref{diag}).  Therefore, by the
canonical bundle formula (\cite[p.~89]{Laz04II}) over
$X=\P(V)\stackrel{\pi}{\>}Y$,
$$-K_X=\sO_{V}(n)\ts\pi^*(\det V^*), $$
where $n$ is the rank of $V$, we have \be -K_X\cdot C&=&\deg_C
\left(\mu^*(-K_X)\right)\\&=&\deg_C\left(\mu^*(\pi^*(\det
V^*))\right)+\deg_C\left(\mu^*\sO_{V}(n)\right)\\
&=&-\deg_C\nu^*(V)+ n\deg_C (L), \ee since $L=\mu^*\sO_V(1)$. Hence,
by the inequality (\ref{sta}), we obtain $-\deg_C\nu^*(V)+ n\deg_C
(L)>0$, which implies $-K_X\cdot C>0$. Therefore,  $-K_X$ is
strictly nef.

$(1)\Longrightarrow (2).$ Let $\nu:C\>Y$ be a finite morphism where
$C$ is an irreducible curve.
  Let $\mu:C\>X$ be the map given  in the configuration (\ref{diag}) and $\tilde C= \mu(C)$.
  Then
\be -K_X\cdot \tilde C&=&\deg_{ C}
\left(\mu^*(-K_X)\right)\\
&=&-\deg_C\nu^*(V)+ n\deg_C (L). \ee
 If $-K_X$ is strictly nef, we
have $-K_X\cdot \tilde C>0$, and so (\ref{sta}) holds. \eproof

\noindent As an application of Lemma \ref{cri}, we obtain

 \bproposition\label{aa}  Let $X=\P(V)$ be a projective bundle of a
nef vector bundle $V$ over a Calabi-Yau manifold $Y$. If there
exists a line bundle $L\in Pic^0(Y)$ such that $H^0(X, V^*\ts L)\neq
0 $, then $-K_X$ is not strictly nef. \eproposition

\bproof Let $\tilde V=V\ts L^*$ and $\tilde X=\P(\tilde V)$. It is
obvious that $X$ is isomorphic to $\tilde X$. Since $L\in Pic^0(Y)$
and $V$ is nef, $\tilde V$ is nef. By using similar arguments as in
Corollary \ref{a}, the nowhere vanishing section $\sigma\in
H^0(Y,\tilde V^*)$ defines a trivial quotient line bundle of
 $\tilde V$. By Lemma \ref{cri}, $-K_{\tilde
X}$ can not be strictly nef. Hence, $-K_X$ is not strictly nef.
\eproof

\subsection{Projective bundles over elliptic curves} In this
subsection, we investigate geometric properties of anti-canonical
bundles of projective bundles over elliptic curves.

\bproposition\label{bb} Suppose $V$ is an indecomposable vector
bundle over an elliptic curve $S$ and $\deg(V)=0$. Then the
anti-canonical line bundle of $\P(V)$ is not strictly nef.
\eproposition \bproof Suppose $V$ has rank $n$. Since $V$ is
indecomposable and $\deg(V)=0$,  by a fundamental result of Atiyah
(\cite[Theorem~5]{Ati57}), there exists a vector bundle $F_n$ such
that $\deg(F_n)=0$, $H^0(S,F_n)\neq 0$ and $V=F_n\ts \det V$. On the
other hand, by Hirzebruch-Riemann-Roch theorem for vector bundles
over elliptic curves, we have \beq \deg(F_n)=\dim
H^0(S,F_n)-H^1(S,F_n).\eeq Hence $H^1(S,F_n)\neq 0$. In particular,
$H^0(S,F_n^*)\neq 0$. Now we consider the variety $X:=\P(F_n)$.  It
is well-known that $\P(V)$ is isomorphic to $\P(F_n)$ since
$V=F_n\ts L$. Suppose $\P(V)$ has strictly nef anti-canonical line
bundle, then so is $-K_X$. By the projection formula again \beq
-K_X=\sO_{F_n}(n)\ts \pi^*(\det F^*_n),\eeq $\sO_{F_n}(1)$ is a
(strictly) nef line bundle, and so $F_n$ is a nef vector bundle. In
summary, $F_n$ is a nef vector bundle over elliptic curve $S$ with
$H^0(S,F_n^*)\neq 0$. However, by Proposition \ref{aa}, $\P(F_n)$
can not have strictly nef anti-canonical bundle. This is a
contradiction. \eproof

\blemma\label{el} Let $X=\P(V)$ be a projective bundle over a smooth
curve $S$. If $-K_X$ is strictly nef, then $S\cong \P^1$. \elemma

\bproof Suppose $S$ is an elliptic curve. We shall use Atiyah's work
\cite{Ati57} on the classification of vector bundles over elliptic
curves to rule out this case. Suppose $V$ has rank $n$.

 \noindent \emph{Case $1$.} Suppose $V$ is indecomposable and
$\deg(V)=0$. This case is ruled out by Proposition \ref{bb}.\\

\noindent \emph{Case $2$}. Suppose $V$ is indecomposable and
$\deg(V)\neq 0$. There exists an \'etale base change $f:Y\>S$ of
degree $k$ where $k$ is an integer such that $n|k$, and $Y$ is also
an elliptic curve. Suppose $X'=\P(f^*V)$, then we have the
commutative diagram
$$\xymatrix{
  X' \ar[d]_{\pi'} \ar[r]^{f'}
                & X\ar[d]^{\pi}  \\
  Y  \ar[r]^{f}
                & S.            }
$$ Hence  $-K_{X'}$ is also strictly nef.
It is obvious that $\rank(f^*V)=\rank(V)=n$, and $\deg(f^*V)=k
\deg(V)$. Let $\ell$ be an integer defined as \beq
\ell=\frac{\deg(f^*V)}{\rank(f^*V)},\eeq
 and $F$ be a line bundle over $Y$ such that $\deg(F)=-\ell$.  Now we set $\tilde V=f^*V\ts F$, then $\deg(\tilde V)=0$.
 A well-known result of Atiyah asserts that an indecomposable vector bundle over an
 elliptic is semi-stable and so $V$ is semi-stable (e.g. \cite[Appendix~A]{Tu93}). Therefore $f^*V$ is
 semi-stable (e.g. \cite[Lemma~6.4.12]{Laz04II}) and so is $\tilde
 V$. Let $\tilde V=\ds \tilde V_i$ where $\tilde V_i$ are
 indecomposable. Since $X'=\P(f^*V)\cong \P(\tilde V)$, by
 projection formula, we have
 $$-K_{X'}=\sO_{\tilde V}(n)\ts \pi^*(\det \tilde V).$$
Since $\deg(\tilde V)=0$, $-K_{X'}$ is strictly nef if and only if
$\sO_{\tilde V}(1)$ is strictly nef. Hence $\tilde V$ is strictly
nef. By Proposition \ref{directsum}, every $\tilde V_i$ is strictly
nef. It is obvious that $\deg(\tilde V_i)=0$. Therefore, $\P(\tilde
V_i)$ is a projective bundle over elliptic curve $C$ with strictly
nef anti-canonical line bundle. By Proposition \ref{bb}, it is
impossible.\\

 \noindent \emph{Case $3$}. In the general case,  suppose $V=\ds_{i} V_i$ where $V_i$ are indecomposable vector bundles
 over $S$. Since $V_i$ is indecomposable,  by Case $1$ and Case $2$, the anticanonical line bundle of $\P(V_i)$ is not strictly nef. By Lemma \ref{cri}, for each $i$,  there exists an irreducible curve $C_i$,
 a finite morphism $\nu_i:C_i\>S$ and a line bundle quotient $L_i$ of $\nu_i^*(V_i)$ such that
 \beq \deg_{C_i} (L_i)\leq \frac{\deg(\nu_i^*(V_i))}{\text{rank}(V_i)}.\eeq
On the other hand, it is an elementary fact that  there exists at
least one
 summand $V_i$ such that
 \beq \frac{\deg(V_i)}{\text{rank}(V_i)}\leq \frac{\deg(V)}{\text{rank}(V)}=\frac{\sum_i\deg(V_i)}{\sum_i \text{rank}(V_i)}.\eeq
For such $V_i$, we have \beq \deg_{C_i} (L_i)\leq
\frac{\deg(\nu_i^*(V_i))}{\text{rank}(V_i)}\leq
\frac{\deg(\nu^*_i(V))}{\rank(V)}.\eeq Note that $L_i$ is also a
line bundle quotient of $\nu_i^*(V)$, hence by Lemma \ref{cri}
again, the
anti-canonical line bundle  $\P(V)$ can not be strictly nef.\\

 Hence, we complete the proof that if $\P(V)\>S$ has
strictly nef anti-canonical bundle, then $S$ is not an elliptic
curve. On the other hand, by a result of Miyaoka \cite{Miy93} (see
also \cite[Corollary~3.14]{Deb01}), if $-K_X$ is nef, then $-K_S$ is
nef. Now we deduce that $S\cong \P^1$.
 \eproof

\btheorem\label{main50}  Let $V$ be a vector bundle over an elliptic
curve $C$. If $V$ is strictly nef, then $V$ is ample. \etheorem

\bproof Since $V$ is also nef, $\deg(V)\geq 0$.  Suppose
$\deg(V)=0$. Let $X=\P(V)$ and $\pi:\P(V)\>C$ be the projection.
Then
$$-K_{X/C}=\sO_V(n)\ts \pi^*(\det V^*).$$
Since $K_C=\sO_C$ and $\deg(V)=0$, $-K_X$ is strictly nef if and
only if $\sO_V(1)$ is strictly nef. By Lemma \ref{el}, $-K_X$ can
not be strictly nef. Hence, $\deg(V)>0$.

 Let $V=\ds V_i$ where $V_i$ are indecomposable. By Proposition
 \ref{directsum}, each $V_i$ is strictly nef and $\deg(V_i)>0$. Since
 indecomposable vector bundles are also semi-stable
 (\cite[Appendix~A]{Tu93}), and a semi-stable vector bundle is ample
 if and only if it has positive degree (\cite[Theorem~6.4.15]{Laz04II}),  we deduce $V_i$ is ample.
 Therefore $V$ is ample.
  \eproof

\vskip 1\baselineskip

\section{The proofs of Theorem \ref{main1},  Theorem \ref{main2} and Theorem \ref{main3}}

\noindent In this section, we will prove Theorem \ref{main1},
Theorem \ref{main2} and Theorem \ref{main3}.

\btheorem\label{main10} Let $X$ be a smooth projective variety. If
$\Lambda^{r}TX$ is strictly nef, then $X$ is uniruled. \etheorem

\bproof Since $\Lambda^{r}TX$ is strictly nef,  $-K_X$ is nef.

  \emph{Case 1}. Suppose there exists an ample divisor $H$ such that $-K_X\cdot
  H^{n-1}>0$, then $X$ is uniruled. Indeed, let $\omega$ be a
  K\"ahler metric in the class $c_1(H)$, then
  \beq \int_X c_1(X)\wedge \omega^{n-1}>0.\eeq
That means the total  scalar curvature of the metric $\omega$ is
strictly positive. Based on the fundamental results in
\cite{BDPP13}, Heier-Wong proved in \cite[Theorem~1]{HW12} that if a
projective variety $X$ has a smooth K\"ahler metric with positive
total scalar curvature, then $X$ is uniruled.

\emph{Case 2}. If $K_X\cdot H^{n-1}=0$ for every ample divisor $H$,
then $K_X$ is numerically trivial (e.g. \cite[3.8, p.~69]{Deb01} ),
i.e. $X$ is a Calabi-Yau manifold. By Yau's theorem (\cite{Yau78}),
the exists a K\"ahler metric $\omega$ such that $Ric(\omega)=0$,
i.e. $T X$ is
 Hermitian-Einstein. Since $\Lambda^{r}TX$ is nef and $c_1(X)=0$, $\Lambda^{r}TX$ is numerically
 flat. (A vector bundle $E$ is called numerically flat (\cite[Definition~1.17]{DPS94}), if both $E$ and $E^*$ are nef, or equivalently, both $E$ and $\det E^*$ are nef.)

If $r=\dim X$, then $\Lambda^{r}TX=-K_X$ is numerically trivial
which is not strictly nef. Hence $1\leq r \leq \dim X-1$. Let $E=TX$
and $\tilde E=\Lambda^r TX$. Then $\tilde E$ is nef and $c_1(\tilde
E)=0$, by \cite[Corollary~2.6]{DPS94} \beq \int_X c^2_1( \tilde
E)\wedge \omega^{n-2}=\int_X c_2(\tilde E)\wedge \omega^{n-2}=0.\eeq
Since $c_1(\tilde E)=0$, $c_1(E)=0$. Let $\lambda_1,\cdots,
\lambda_n$ be the Chern roots of $E$, then $c_2(\tilde E)$ is a
symmetric polynomial in $\lambda_1,\cdots, \lambda_n$ of degree $2$.
In particular, \beq c_2(\tilde E)= a c_1^2(E)+b c_2(E)\eeq since
$c_1(E)$ and $c_2(E)$ are elementary symmetric polynomials. Note
that $a$ and $b$ are constants depending only on $r$ and the rank of
$n$ of $E$. A simple computation shows $b\neq 0$. Hence, \beq \int_X
c_2(E)\wedge\omega^{n-2}=0.\eeq Hence, $X$ is a Calabi-Yau manifold
with \beq \int_X c^2_1( X)\wedge \omega^{n-2}=\int_X c_2(
X)\wedge\omega^{n-2}=0,\eeq By Yau's criterion
\cite[Theorem~1.4]{Yau77} (see also \cite[Corollary~9.6]{Zhe00}),
$TX$ is Hermitian flat and up to a finite \'etale
 cover, $X$ is isomorphic to a torus. This is a contradiction since $\Lambda^r
 TX$ is strictly nef.
\eproof

\btheorem Let $X$ be a smooth projective variety with dimension $n$.
If $TX$ is  strictly nef, then $X\cong \P^n$. \etheorem

\bproof Since $TX$ is  strictly nef, $TX$ is nef. By the structure
theorem (\cite[Main Theorem]{DPS94}) of projective manifolds with
nef tangent bundles, there exists a finite \'etale cover $\pi:\tilde
X\>X$ such that the Albanese map $\alpha:\tilde X\>A$ is a smooth
morphism and the fibers are Fano manifolds with nef tangent bundles.
Note that by Theorem \ref{main1}, $\dim A<\dim X$. Since $\pi$ is a
finite morphism, $T\tilde X$ is also strictly nef. Now consider the
exact sequence \beq 0\>T_{\tilde X/A}\>T{\tilde
X}\>\alpha^*(TA)\>0.\eeq Since $T\tilde X$ is strictly nef, it can
not have a trivial quotient. Hence $\dim A=0$, and  $X$ is Fano. For
an arbitrary rational curve $\nu:\P^1\>X$, \beq \nu^*TX=\bds_{i=1}^n
\sO_{\P^1}(a_i), \eeq where $a_i\geq 1$ for each $i$. Therefore
$\deg(\nu^*K^{-1}_X)\geq n.$ If $\deg(\nu^*K_X^{-1})=n$, then
$a_i=1$ for all $i$ which can imply
 the tangent morphism $T\P^1\>\nu^*TX$ is zero, which is a
 contradiction. That means $\deg(\nu^*K^{-1}_X)\geq n+1$. Recall that the
pseudo-index of a Fano variety $X$ is defined as:
$$i_X=\min\{-K_X\cdot C|\text{ C is a rational curve} \},$$
   and so the pseudo-index $i_X\geq n+1$. Hence,
 $X\cong \P^n$.
 \eproof

\noindent

\noindent The following result will be used frequently in the proof
of Theorem \ref{main3}.

\blemma\label{number} Let $X$ be a smooth projective variety with
dimension $n\geq 3$. Suppose $\Lambda^2T_X$ is  strictly nef. \bd
\item For any rational curve $C$,   we have \be -K_X\cdot C\geq
n.\ee
\item If $T_X$ is not nef along a rational curve $C$, then
\be -K_X\cdot C\geq 2n-3.\ee \ed \elemma

\bproof $(1)$ Over the rational curve $C$, $T_X$ splits as \beq
 T_X|_{ C}=\bigoplus_{b_j\leq0}\mathcal{O}(b_j)\bigoplus\left(\bigoplus_{a_i>0}\mathcal{O}(a_i)\right),\eeq
Since $\Lambda^2 T_X$ is  strictly nef, by Lemma \ref{cri0} \beq
T_X|_{C}=\mathcal{O}(-b)\bigoplus\left(\bigoplus_{i=1}^{n-1}\mathcal{O}(a_i)\right),\eeq
 where $0<a_1\leq a_2 \leq \cdots \leq a_{n-1}$  and $a_1-b>0$.

 If $b>0$, then $a_1\geq2$ and \beq -K_X\cdot C=(a_1-b)+a_2+\cdots +a_{n-1}\geq
 1+2(n-2)=2n-3\geq
 n\eeq since $n\geq 3$.

 If $b=0$, since from the tangent morphism $0\>\sO(2)\>TX|_C$ of
 $C\>X$, there exists at least one $a_i\geq 2$
and so $-K_X\cdot C\geq n$.

$(2)$. If $T_X$ is not nef along a rational curve $C$, then we have
$b>0$ and so $-K_X\cdot C\geq 2n-3$. \eproof

\noindent One of the key ingredients in the proof of  Theorem
\ref{main3} is

\btheorem\label{keytheorem} Let $X$ be a smooth projective variety
of dimension $n> 4$. Suppose that $\Lambda^2 T_{X}$ is strictly nef,
then $c_1^n(X)\neq 0$. \etheorem

\noindent This is the most complicated issue in the proof of Theorem
\ref{main3}, since as a priori, it might be possible that
$c_1^n(X)=0$ (c.f. the Hermitian flat strictly nef vector bundle in
Example \ref{example}). Based on Theorem \ref{keytheorem}, we are
ready to prove Theorem \ref{main3}.

\btheorem \label{wedge} Let $X$ be a smooth projective variety of
dimension $n> 4$. Suppose that $\Lambda^2 T_{X}$ is strictly nef,
then $X$ is isomorphic to $\mathbb{P}^n$ or a quadric $\Q^n.$
\etheorem

\bproof  Since $\det (\Lambda^2 T_{X})=(n-1)(-K_X),$ the strictly
nefness of $\Lambda^2 T_{X}$ implies that $-K_X$ is nef and
$c^n_1(X)\geq 0$. By Theorem \ref{keytheorem}, $-K_X$ is big, i.e.
$c^n_1(-K_X)>0$. By Kawamata-Reid-Shokurov base point free theorem,
$-K_X$ is semiample.
 Then there exists $m$ big enough such that $\phi=|-m K_X|$ is a morphism. Since
$-K_X$ is big, there is a positive integer $\tilde m$ such that
$$\tilde m(-K_X)=D+L$$
where $D$ is an effective divisor and $L$ is an ample line bundle.
If $-K_X$ is not ample, there exists exists a curve $C$ contracted
by $\varphi$ which means $-K_X\cdot C=0$. Therefore,
$$ D\cdot C=-L\cdot C<0.$$ Let $\Delta=\eps D$ for some small $\eps>0$, then
$(X,\Delta)$ is a klt pair and $K_X+\Delta$ is not $\varphi$-nef.
Then by the relative Cone theorem(e.g. \cite[Theorem~3.25]{KM98})
for log pairs, there exists a rational curve $\tilde{C}$ contracted
by  the morphism $\phi$ . By Lemma \ref{number} \beq
T_X|_{\tilde{C}}=\mathcal{O}(-b)\bigoplus\left(\bigoplus_{i=1}^{n-1}\mathcal{O}(a_i)\right),\eeq
where $a_1-b\geq1.$ Hence, \begin{eqnarray}
       (-K_X)\cdot \tilde C\geq 1
     \end{eqnarray}
 which is a
contradiction since $\tilde{C}$ is contracted by $\varphi.$
Therefore,  $-K_X$ is ample and $X$ Fano. By the computation in
Lemma \ref{number}, we know the pseudo-index $i_X\geq n$. Thanks to
\cite{Miy04} or \cite[Theorem~C]{DH}, $X$ is isomorphic to
$\mathbb{P}^n$ or a quadric $\Q^n.$\eproof

\noindent In the rest of this section, we shall prove  Theorem
\ref{keytheorem}.   We need several technical lemmas to complete the
proof.

\blemma \label{bend and break} Let $X$ be a smooth projective
variety of dimension $n\geq 4$. Suppose  $\Lambda^2 T_X$  is
strictly nef. Let $Y$ be a prime divisor of $X$  and
$\mathcal{L}=\mathcal{O}_X(Y)$.  Suppose        $ Y $ is rationally
chain connected and     $ \mathcal{L}|_{Y} $ is numerically trivial.
Then $Y$ is a smooth variety.

\elemma \bproof We will follow some ideas of \cite[Lemma~4.12]{SW}
and use the techniques of \cite[Theorem~1.3, Chapter~II]{Kol01}. Let
$\sJ$ be the  ideal sheaf of $Y$. Then
$\mathcal{J}/\mathcal{J}^2\simeq \mathcal{L}^{-1}|_{Y}.$   There
exists an exact sequence: \beq
   \mathcal{L}^{-1}|_{Y}\>\Omega_{X}|_{Y} \>\Omega_{Y} \> 0. \label{bk}\eeq
Let $U$ be the smooth locus of $Y$. Suppose that there is a rational
curve $f: \mathbb{P}^1\longrightarrow Y$ whose image has nonempty
intersection with $U$.\\

  \emph{Claim} $1$: If $f^*(T_X)$ is nef, then $Y$ is smooth along $f(\mathbb{P}^1),$ i.e. $f(\mathbb{P}^1)\subseteq
  U$.\\

    Actually, by restricting the exact sequence (\ref{bk}) to $f(\P^1),$ we have the following right exact sequence:
  \beq
     \mathcal{O}\stackrel{d}{ \>} f^*(\Omega_X) \>  f^*(\Omega_Y) \> 0.
    \label{bk2}\eeq
   Note that $d$ is generically split, which means $d$ is injective on the  dense subset $f(\mathbb{P}^1)\bigcap U\neq \emptyset$,  so $\Ker(d)$ is a torsion free (hence locally free) subsheaf of generic rank $0,$ then $\Ker(d)=0$ and $d$ is injective.
     So the right exact sequence (\ref{bk2}) is also left exact. Thus, the morphism $d$ defines a global section $s$
     of $H^0(\P^1, f^*(\Omega_X))$. On the other hand, since $f^*(T_X)$ is nef, by \cite[Proposition~1.16]{DPS94},
      $s$ is nowhere vanishing. Then $ f^*(\Omega_Y)$ is of constant rank of $n-1$.  So by \cite[Theorem~8.15, Chapter~II]{Har}, $Y$ is  smooth along $f(\mathbb{P}^1).$\\

    \emph{Claim} $2$: If $f^*(T_X)$ is not nef, then $\dim_{[f]} Hom(\mathbb{P}^1, Y|_{\{0, \infty\}})\geq 2.$\\

   Since $f^*(\Lambda^2T_X)$ is  strictly nef and $f^*(T_X)$ is not nef,
   $$f^*(T_X)=\mathcal{O}(-b)\bigoplus\left(\bigoplus_{i=1}^{n-1}\mathcal{O}(a_i)\right),$$ where $b\geq1 $ and $a_i\geq 2.$ By the exact sequence (\ref{bk2}) and the fact
$f^*(\mathcal{L}^{-1})=\mathcal{O}$,  we obtain the long exact
sequence:
  \be
    0 \> Hom( f^*(\Omega_Y), \mathcal{O}) &&\> Hom( f^*(\Omega_X), \mathcal{O}) \> Hom(\mathcal{O} , \mathcal{O})\\&& \> \Ext^1( f^*(\Omega_Y), \mathcal{O})
     \>\Ext^1( f^*(\Omega_X), \mathcal{O})\> \Ext^1(\mathcal{O} ,
     \mathcal{O}).\ee
    Then by \cite[Theorem~2.16, Chapter~I]{Kol01}, \be
    \dim_{[f]}Hom(\mathbb{P}^1,Y) &\geq& \dim Hom (f^*(\Omega_Y), \mathcal{O})-\dim \Ext^1( f^*(\Omega_Y), \mathcal{O})\\
     &=&\chi(\mathbb{P}^1, f^*(T_X))-\chi(\mathbb{P}^1,\mathcal{O}) \\ &=&(1-b)+\sum^{n-1}_{i=1}(1+a_i)-1\\
     &\geq&
     3n-4.\ee
    Since $Hom(\mathbb{P}^1, Y|_{\{0, \infty\}})$ is the fibre of the restricting morphism
     \beq Hom(\mathbb{P}^1, Y)\longrightarrow Hom(\{0,\infty\}, Y)\simeq Y\times
     Y,\eeq
   we have  \beq \dim_{[f]}Hom(\mathbb{P}^1, Y|_{\{0, \infty\}})\geq \dim_{[f]} Hom(\mathbb{P}^1, Y)-2\dim Y\geq
   n-2\geq2.\eeq
\vskip1\baselineskip

   \emph{Claim} $3$: If there is a connected rational chain $\Gamma$ such that  $\Gamma\bigcap U\neq \emptyset$, then $\Gamma\subseteq U$. \\

   We prove it by induction on the degree $-K_X\cdot \Gamma$.
   If $-K_X\cdot \Gamma=n$, by Lemma \ref{number},  $\Gamma$ is irreducible and $T_X$ is nef along $\Gamma$ since the dimension $n\geq 4$ and $n<2n-3$.
    By Claim $1$, we know that $\Gamma\subseteq U.$\\

  By induction,  we assume that any connected rational chain $\Gamma$ which has nonempty intersection with $U$ and  $-K_X \cdot \Gamma<N$ is contained in $U$.
  We want to show that if there is a connected rational chain $\Gamma$ such that  $\Gamma\bigcap U\neq \emptyset$ and $-K_X\cdot \Gamma= N$, then $\Gamma\subseteq U$.   \\

  If $\Gamma$ is reducible, i.e. $\Gamma=\Gamma_1+\Gamma_2,$ since $\Gamma\bigcap U\neq\emptyset,$ there is at least one component which has nonempty intersection with
  $U$,
  namely $\Gamma_1.$ Then $-K_X\cdot \Gamma_1< N,$ by induction, $\Gamma_1\subseteq U.$
  Since $\Gamma$ is connected, $\emptyset\neq\Gamma_2\bigcap\Gamma_1\subseteq\Gamma_2\bigcap U.$ So by induction again, $\Gamma_2\subseteq U.$\\

  If $\Gamma$ is irreducible, by Claim $1,$ we may assume  $T_X$ is not nef along $\Gamma.$ Then fix two points $x$ and $y$ of
  $\Gamma$,
  where $x\in \Gamma\bigcap U.$ Then by Claim $2,$ $\dim_{[f]} Hom(\mathbb{P}^1,Y|_{\{0,\infty\}})\geq 2$.
   Therefore, by \cite[Theorem~4]{Mor79} or \cite[Proposition~3.2]{Deb01},
    $\Gamma$ breaks into a connected rational chain with $x$ and $y$ fixed. We denote it by $\Gamma\equiv\Gamma_1+\Gamma_2.$
    By the same argument as above, $\Gamma_1, \Gamma_2 \subseteq U,$ hence $y\in U.$ As $y$ moves along $\Gamma,$ $\Gamma\subseteq U.$\\

 \noindent  Then Lemma \ref{bend and break} follows directly from Claim $3.$\eproof

\blemma \label{smooth fibration} Let $X$ be a smooth projective
variety of dimension $n>4$. Suppose  $\Lambda^2 T_X$  is strictly
nef. If $T_{X}$ is not nef, then  there exists a Mori elementary
contraction $\beta:X\longrightarrow B$ such that $\beta$ is a
projective bundle. \elemma

\bproof  By Theorem \ref{main10} and Lemma \ref{number}, there
exists rational curve $\sigma:C\>X$ such that  $\upsilon^*(-K_X)$ is
ample over $C$. Then by the ``Cone theorem'' (e.g.
\cite[Theorem~3.15]{KM98}), there exist countably many rational
curves $\Gamma_i$ such that
 $$ \overline{NE(X)}=\overline{NE(X)}_{K_X\geq0}+\sum_{i}\mathbb{R}^+ [\Gamma_i]$$
 where $0<-K_X\cdot \Gamma_i\leq n+1.$
 Since $\Lambda^2 T_X$ is  strictly nef,
 $$T_X|_{\Gamma_i}=\mathcal{O}(-b)\bigoplus\left(\bigoplus_{i=1}^{n-1}\mathcal{O}(a_i)\right),$$
 Then by Lemma \ref{number}, $\Gamma_i$ is a free curve, i.e. $b\leq 0$, otherwise $-K_X\cdot {\Gamma_i}\ge 2n-3>n+1$.
 Since deformations of a free curve cover a dense subset of $X$ {and $Locus(\Gamma_i)$  is closed,  $Locus(\Gamma_i)=X$. So for an arbitrary $\Gamma_i$,
  the associated Mori elementary contraction is of fibre type}.
  We pick one of these extremal curves, namely $\Gamma$ and denote the associated contraction by $$\beta:X\longrightarrow B.$$
For an arbitrary smooth fibre $F,$ there exists an exact sequence$$
  0\> T_F \>T_{X}|_{F} \>\mathcal{O}_F^{\oplus\dim B} \> 0 .$$
Since $\Lambda^2T_X$ is strictly nef, we deduce $\dim(B)\leq 1$.
Actually, $B$ is not a point, otherwise the Picard number
$\rho(X)=1$ and $X$ is Fano. By Lemma \ref{number},  $X$ is
$\mathbb{P}^n$ or $\Q^n$. So $B$ is a curve and
$F\simeq\mathbb{P}^{n-1}.$ Now we aim to show that $T_X$ is
$\beta$-nef, then by \cite[Theorem ~4.2]{SW}, $\beta$ is a smooth
fibration. \\

\emph{Claim} $1:$ Let $F'_0$ be a fibre (possibly singular) of $\beta:X\>B$ and $F_0=red(F'_0),$  then $F_0$ is smooth and irreducible.\\

Since  fibres  of $\beta$ are equi-dimensional, by
\cite[Theorem~3.5.3, Chapter ~4]{Kol01}, every fibre is rationally
chain connected. Now suppose that $F_0$ is reducible and $F_0=
F_1+F_2.$ Note that for an arbitrary curve $C$ in $F_0,$ by the
``Cone Theorem'', $C$ is numerical to some multiple of $\Gamma,$
which means any fibre curve could be deformed into smooth fibres.
Then for  an arbitrary rational curve $C_1$ contained in $F_1$ which
has non-empty intersection with $F_2,$  $C_1\cdot F_2=0,$ then $C_1$
is contained in $F_2.$ So $F_0$ is
irreducible. It  is obvious  that $F_0$ satisfies the assumptions of Lemma \ref{bend and break}, so $F_0$ is smooth.\\

\emph{Claim} $2:$ $T_{X}$ is $\beta$-nef.\\

Since $F_0$ is smooth, for an arbitrary curve $g:
\tilde{C}\longrightarrow F_0,$ by the same arguments as in Lemma
\ref{bend and break}
 \beq
 0\>g^*(\mathcal{J}/\mathcal{J}^2) \> g^*( \Omega_X) \> g^*(\Omega_{F_0})\> 0
 \eeq
 is exact. Note that $\mathcal{J}/\mathcal{J}^2\simeq \mathcal{O}_X(F_0)^{-1}|_{F_0}$ and $\mathcal{O}_X(F_0)|_{F_0}$ is a numerically trivial line bundle.
 By the induced exact sequence
\beq 0\>g^*\left(\Lambda^2T_{F_0}\right) \> g^*\left(\Lambda^2
T_X\right) \> g^*(T_{F_0})\otimes g^*(\mathcal{O}_X(F_0)|_{F_0})\>
0, \eeq
we obtain $g^*(T_{F_0})$ is nef. Then as an extension of $g^*(T_{F_0})$ with a numerically trivial line bundle, $g^*(T_X)$ is nef.\\

 \noindent In conclusion, Lemma \ref{smooth fibration} is a direct corollary of Claim
 $2$.
\eproof

\btheorem\label{nef}  Let $X$ be a smooth projective variety of
dimension $n>4$. If  $\Lambda^2 T_X$  is  strictly nef, then $T_X$
is nef. \etheorem \bproof
 Suppose $T_X$ is not nef, then by
Lemma \ref{smooth fibration}, $X$ is a projective bundle over some
smooth curve $B$. Since $-K_X$ is nef,  $B$ is an elliptic curve or
$\mathbb{P}^1.$

If $B$ is an elliptic curve, then by the relative tangent bundle
exact sequence \beq 0 \> T_{X/B} \> T_{X} \>\beta^*(T_B) \> 0 ,\eeq
and the induced exact sequence \beq 0 \> \Lambda^2 T_{X/B} \>
\Lambda^2T_{X} \>T_{X/B} \> 0 ,\label{relatvie}\eeq we deduce that
$T_X$ is nef since $\Lambda^2TX$ and $T_{X/B}$ are nef. This is a
contradiction.

Suppose $B=\P^1$. Let $X=\mathbb{P}(V)$ for some vector bundle $V,$
we may assume that $V=\bigoplus_{i=1}^n\mathcal{O}(a_i)$ where
$0<a_1\le a_2\le\cdots\le a_n.$ Let $\mathcal{L}=\mathcal{O}(a_1),$
then the quotient morphism
$$\bigoplus_{i=1}^n\mathcal{O}(a_i)\longrightarrow\mathcal{O}(a_1)\longrightarrow
0$$ defines a section $\sigma: B\longrightarrow X.$ Indeed, let
$\nu:B\>B$ be the identity map.  By (\ref{diag}), this quotient
determines a section of $\beta:\P(V)\>B$, namely, a  morphism
$\sigma: B\longrightarrow \tilde X$ such that $\beta\circ
\sigma=\nu=Id$ with $\sL=\sigma^*(\sO_V(1))$.

Since $\beta: X\longrightarrow B$ is a projective bundle,  one has
the Euler sequence \beq
  0 \> \sO_X \> \beta^*(V^*)\ts \sO_V(1) \> T_{X/B} \>  0
 \label{Euler} \eeq
 and
$$\sigma^*\left(\beta^*(V^*)\ts \sO_X(1)\right)=\sigma^*\circ\beta^*(V^*)\otimes\mathcal{L}=\bigoplus_{i=1}^n\mathcal{O}(a_1-a_i).$$
 Let $\tilde{a_i}=a_1-a_i,$ then $\tilde{a_i}\le 0$ and $0=\tilde{a_1}\ge \tilde{a_2}\ge\cdots\ge\tilde{a_n}$. Let $\sigma^*(T_{X/B})=\bigoplus_{j=1}^{n-1}\mathcal{O}(b_j),$
 where $b_1\le\cdots\le b_{n-1}$.  By pulling-back the exact sequence (\ref{Euler}),  we obtain the following exact sequence:
 \beq
   0 \>\sO_B \> \bigoplus_{i=1}^n\mathcal{O}(\tilde{a_i})\>\bigoplus_{j=1}^{n-1}\mathcal{O}(b_j)
   \>0.\label{Euler2}
   \eeq
By performing square wedge product of the exact sequence
(\ref{Euler2}), we obtain another exact sequence: \beq
  0 \> \bigoplus_{j=1}^{n-1}\mathcal{O}(b_j)\stackrel{d}{ \>}\bigoplus_{i<j}\mathcal{O}(\tilde{a_i}+\tilde{a_j}) \>\bigoplus_{i<j}\mathcal{O}(b_i+b_j) \> 0
  \eeq
 By the injectivity of $d$,  we know $b_j\le \tilde{a_1}+\tilde{a_2}\le 0$ for arbitrary $j.$ Since the sequence $$0 \> T_{X/B} \> T_{X} \>\beta^*(T_B) \> 0$$
  is exact, $\deg \sigma^*(-K_X)=\Sigma b_j +2\le 2.$ Since $\Lambda^2T_X$ is  strictly nef, we have already obtained the inequality $\deg\sigma^*(-K_X)\ge n$ in Lemma \ref{number}.
   This is a contradiction.\eproof

\noindent Now we complete the proof of Theorem \ref{keytheorem}, by
using Theorem \ref{main2}, Theorem \ref{main5} and the lemmas
described above.\\

\noindent \emph{The proof of Theorem \ref{keytheorem}.}
 By Theorem \ref{nef}, $TX$ is nef. Suppose $c_1(X)^n=0$.  By
Theorem \cite[Main Theorem]{DPS94}, $X$  has a finite \'etale cover
$\tilde X$ with irregularity $q(\tilde X )>0$ and the Albanese
morphism $\alpha: \tilde X\longrightarrow Alb(\tilde X)$ is a smooth
fibration. Moreover, the fibers $F$ are Fano manifolds with nef
tangent bundles. Let $C:=Alb(\tilde X)$. Since $T_{\tilde X/C}|_F=
T_F$ and $\alpha^*(T_C)|_{F}$ is a trivial vector bundle, we obtain
the following exact sequence from the relative exact sequence of
tangent bundles: \beq
  0 \> T_F \> T_{\tilde X}|_{F} \> \mathcal{O}_F^{\oplus\dim C} \> 0.
  \eeq
Since $\Lambda^2T_{X}$ is  strictly nef, $\Lambda^2T_{\tilde X}$ is
strictly nef.  By Lemma \ref{cri0}, we know $\dim C=1$. From  the
induced exact sequence \beq
  0 \> \Lambda^2 T_F \> \Lambda^2 T_{\tilde X}|_{F} \> TF \> 0.
  \eeq
we know  $T_{F}$ is  strictly nef. By Theorem \ref{main2}, $F$ is
isomorphic to $\mathbb{P}^{n-1}.$  Then $\tilde X$ is a
$\P^{n-1}$-bundle over an elliptic curve $C.$  Therefore, $\tilde X$
is a projective bundle over $C$.  Let $\tilde X=\P(V)\>C$.  We have
$ -K_{\tilde X}=\sO_V(n)\ts \pi^*\det (V^*).$  Since $ T\tilde X$ is
nef, and $\P(V)=\P(V\ts L)$ for any line bundle $L$, after a finite
\'etale cover (e.g. Case $2$ of Lemma \ref{el}), we can assume
$\deg(V)=0$ and so $V$ is numerically flat.

By Lemma \ref{flat}, a numerically flat vector bundle  $V$ over an
elliptic curve admits a line bundle quotient $V\> L\> 0,$ where $L$
is  of degree $0$. Let $\nu:C\>C$ be the identity map. By the
configuration in (\ref{diag}), this quotient determines a section of
$\pi:\P(V)\>C$, namely, a  morphism $\sigma: C\longrightarrow \tilde
X$ such that $\pi\circ \sigma=\nu=Id$. By pulling back the relative
tangent bundle sequence as in (\ref{relatvie}), we obtain the
 exact sequence:
\beq 0 \> \sigma^*\left(\Lambda^2T_{\tilde X/C}\right) \>
\sigma^*\left(\Lambda^2T_{\tilde X}\right) \>
\sigma^*\left(T_{\tilde X/C}\right) \> 0 .\label{relatvie2}\eeq
   Since $\sigma^*(\Lambda^2 T_{\tilde X})$ is strictly nef, as its quotient, $\sigma^*(T_{\tilde X/C})$ is strictly nef over $C$.
 As $C$ is an elliptic curve, by Theorem \ref{main50}, $\sigma^*(T_{\tilde X/C})$ is ample.
 Then by pulling back the Euler sequence (\ref{Euler}), we have
$$
    0 \>\mathcal{O}_C  \> \sigma^*(V^*)\otimes L  \> \sigma^*(T_{\tilde X/C})\> 0,$$
since $L=\sigma^*(\sO_V(1))$ as in (\ref{diag}). Since
$\sigma^*(T_{\tilde X/C})$ is ample, $$\deg\left(
\sigma^*(V^*)\otimes L\right)=\deg( \sigma^*(T_{\tilde X/C}))>0.$$
     Since $ V$ is numerical flat and $\deg(L)=0$, we obtain a
     contradiction. The proof of Theorem \ref{keytheorem} is completed.\qed

\blemma\label{flat} Let $C$ be an elliptic curve. \bd \item Let $V$
be a Hermitian flat vector bundle. Then $V=L_1\ds\cdots\ds L_r$
where $L_i\in Pic^0(C)$ and $r$ is the rank of $V$.
\item Let $V$ be a numerically flat vector bundle.  $V$ admits a line bundle quotient $V\> L\> 0,$
where $L\in Pic^0(C)$.

\ed

\elemma

\bproof $(1)$.  By Harder-Narasimhan Theorem (e.g.
\cite[Theorem~2.5, Chapter~V]{Kob87}), if $V$ is Hermitian flat over
an elliptic curve, it is Hermitian-Einstein and so $ V=\ds V_i$
where $V_i$ are stable bundles of degree zero. It is well-known
that, $V_i$ must be of rank one (e.g. \cite[Appendix ~A]{Tu93}),
i.e. $V_i\in Pic^0(C)$.

$(2).$ If $V$ is numerically flat,  $V^*$ is also numerically flat.
By \cite[Theorem~1.18]{DPS94}, a holomorphic vector bundle $V^*$ is
numerically flat if and only if $V^*$ admits a filtration \beq
\{0\}=V_0\subset V_1\subset \cdots \subset V_p=V^*\eeq by vector
subbundles such that the quotients $V_k/V_{k-1}$ are Hermitian flat.
In particular, $V_1$ is Hermitian flat and by $(1)$, $V_1$ has a
degree zero line subbundle $L^*$. Hence $L^*$ is a line subbundle of
$V^*$, i.e. $V\>L\>0$ is a line bundle quotient. \eproof

\end{document}